\documentclass[11pt]{article}
\usepackage{times}
\usepackage{amsmath,amssymb}
\usepackage{theoremref}
\usepackage[utf8]{inputenc}
\usepackage{indentfirst}
\usepackage{graphicx,tikz}
\usepackage{xcolor}
\usepackage{xskak}
\usepackage[T1]{fontenc}
\usepackage{amsthm}
\newtheorem{problem}{Question}

\newtheorem{theorem}[problem]{Theorem}

\newtheorem{conjecture}[problem]{Conjecture}

\title{A Correlation Inequality on Three Functions}
\author{Kada K Williams}

\begin{document}

\maketitle

\begin{abstract}
Let $X$ and $Y$ be upward closed set systems in the lattice of $\{0,1\}^n$. The celebrated Harris-Kleitman inequality implies that if $|X|=\alpha 2^n$, $|Y|=\beta 2^n$, the density of the set of points in exactly one of $X$ and $Y$ is maximal when $X$ and $Y$ are independent, meaning $|X\cap Y|=\alpha\beta 2^n$. Is the same true of three upward closed systems, $X$, $Y$, and $Z$? Suppose $|X|=|Y|=|Z|$. Kahn asked whether the set of points in exactly one of $X$, $Y$, $Z$ has density at most $\frac49$. We answer this question in the negative.
\end{abstract}


\section{Introduction}

Let $(P,<)$ be a finite poset. We say that $X\subset P$ is an upward closed system if whenever $a<b$ and $a\in X$, also $b\in X$. Notice that if $X$ and $Y$ are upward closed, then so is $X\cap Y$.

Of particular interest are upward closed systems in the poset $Q_n$. Our notation is standard, with $[n]=\{1,2,\dots,n\}$ and $Q_n$ is the power set of $[n]$ ordered by $\subset$.

\begin{theorem}[Harris-Kleitman] \label{HKI}
Let $X$ and $Y$ be upward closed set systems on a ground set of $n$ elements. Then the following inequality holds \cite{Har}, \cite{Kle}:
$$\frac{|X\cap Y|}{2^n}\ge \frac{|X|}{2^n}\cdot \frac{|Y|}{2^n}.$$
\end{theorem}

Observe that if $X$ is upward closed, its complement $X^c$ is downward closed, meaning that it is upward closed in the reverse of $Q_n$. We deduce that if $|X|=\alpha 2^n$, $|Y|=\beta 2^n$, then $|X\cap Y|\ge \alpha\beta\cdot 2^n$ and $|X^c\cap Y^c|\ge (1-\alpha)(1-\beta)\cdot 2^n$. Therefore,
$$|X\cap Y^c|+|X^c\cap Y|=2^n-|X\cap Y|-|X^c\cap Y^c|\le (\alpha+\beta-2\alpha\beta)\cdot 2^n.$$

Let $X$, $Y$, and $Z$ be upward closed systems in $Q_n$. Applying the Harris-Kleitman inequality to $X\cap Y$ and $Z$ gives a lower bound on the set of points in all three sets, and similarly for their complements. It is therefore natural to explore the set of points in exactly one of the three systems,
$$S_1=(X\cap Y^c\cap Z^c)\sqcup (X^c\cap Y\cap Z^c)\sqcup (X^c\cap Y^c\cap Z).$$
Let $|X|=\alpha 2^n$, $|Y|=\beta 2^n$, and $|Z|=\gamma 2^n$. Supposing $Q_{n}\cong Q_a\times Q_b\times Q_c$, where $a+b+c=n$, let a set be in $X$ depending only on its $Q_a$ component, and similarly for $Y$ with $Q_b$ and $Z$ with $Q_c$. Then
$$|X\cap Y^c\cap Z^c|=\alpha(1-\beta)(1-\gamma)\cdot 2^n,$$
and similarly for $X^c\cap Y\cap Z^c$ and $X^c\cap Y^c\cap Z$. Specially, if $\alpha=\beta=\gamma=\rho$, then
$$\frac{|S_1|}{2^n}= 3\rho(1-\rho)^2\le \frac49.$$

\begin{conjecture}\label{equal}
Let $X,Y,Z\subset Q_n$ be upward closed set systems of density $\rho$. Then 
$$\frac{|S_1|}{2^n}\le 3\rho(1-\rho)^2.$$
\end{conjecture}

While this conjecture is false for $n=5$, an upper bound is relatively easy to glean for equal-sized systems.

\begin{theorem}\label{equaltrue}
Let $|X|=|Y|=|Z|=\rho\cdot 2^n$. Then $\frac{|S_1|}{2^n}\le 3\rho\cdot \frac{1-\rho}{1+\rho}$.
\end{theorem}

The maximum of $\frac{3\rho (1-\rho)}{1+\rho}$ is $\approx 0.515$, achieved when $\rho=\sqrt2-1\approx 0.414$. This yields a heuristic that if $\rho=\frac38\approx 0.375$, the size of $S_1$ can be nearly maximal.

\begin{conjecture}[Kahn]\label{Kahn}
If $|X|=|Y|=|Z|$, then $\frac{|S_1|}{2^n}\le \frac49$.
\end{conjecture}

This conjecture is due to Jeff Kahn \cite{Nar}. We provide a counterexample in $n=21$ dimensions, where the density of $S_1$ exceeds $0.447$.

\section{Proof that Conjecture \ref{equal} is false}

In $Q_5$, let $X$ and $Y$ contain the sets that have $1$ and $2$ as an element, respectively. To begin with, let $W$ contain the sets with at least $3$ elements in $[5]$. Clearly, $X$, $Y$, and $W$ are upward closed with density $\frac12$. Moreover, $|X\cap Y^c\cap W^c|=|X^c\cap Y\cap W^c|=4$, counting a singleton and three pairs.

The only set in $W\setminus X\setminus Y$ is $\{3,4,5\}$. This can be improved by introducing the sets $\{3,4\}$ and $\{3,5\}$ into $W$, removing the sets $\{1,4,5\}$ and $\{2,4,5\}$ in return. The resulting set system $Z$ is upward closed, and while $X^c\cap Y^c\cap Z$ increased by $2$, there is also a new set contained exactly in $X$ and exactly in $Y$. Thus, $|S_1|=5+5+3=13$, which is larger than $\frac49\cdot 2^5=12$.

\section{The upper bound of Theorem \ref{equaltrue}}

For each $i=0,1,2,3$, let $s_i\cdot 2^n$ be the number of sets in $Q_n$ that occur in $i$ many of $\{X,Y,Z\}$. We normalised such that
$$s_0+s_1+s_2+s_3=1.$$
Theorem \ref{HKI} implies that $|X\cap Y\cap Z|\ge \rho |X\cap Y|$, and similarly for $Y\cap Z$ and $Z\cap X$. Adding these up yields $3(1-\rho)s_3\ge \rho s_2$, where $s_2\ge 0$. The same argument for the complements informs us that $3\rho s_0\ge (1-\rho)s_1$, where $s_1\ge 0$. Finally, the constraint on the size of $X$, $Y$, and $Z$ implies that
$$s_1+2s_2+3s_3=3\rho.$$
Given these constrains on $(s_0,s_1,s_2,s_3)$, what is the maximum of $s_1$?

The strategy for solving such an optimisation problem is the following. As long as the linear equations for $(s_0,s_1,s_2,s_3)$ are underdetermined, it is possible to move this vector along a line while increasing $s_1$, until one of the inequalities becomes an equality. Iterating this observation, we deduce that $s_2$ is maximal with equalities $3s_0=\frac{1-p}{p}s_1$ and $s_2=0$. Solving the linear system, we obtain $s_1\le \frac{3(1-\rho)\rho}{1+\rho}$. The maximum of this bound is $9-6\sqrt 2\approx 0.515$, achieved when $\rho=\sqrt2-1\approx 0.414$.

Can equality hold in the poset $Q_n$? Consider the weighted poset
$$P=\{a,A,p_1,p_2,p_3\},$$
where $a<p_i$ and $p_i<A$, with weights $\frac{(1-p)^2}{1+p}$, $\frac{(1-p)p}{1+p}$, and $\frac{2p^2}{1+p}$ for $a$, $p_i$, and respectively. Theorem \ref{HKI} holds for upward closed systems in $P$, the only remarkable instance being of two systems of size three, where $\left(\frac{2p}{1+p}\right)^2\le p$. It would seem that $P$ cannot be simulated in $Q_n$, but we did not find any obvious reason for this.

\section{A counterexample for Conjecture \ref{Kahn}}

Like Sahi \cite{Sah}, we consider the weighted hypercube $Q_n(p)$, which is the measure space of $n$ independent Bernoulli bits $B(p)$. That means that a set of $k$ elements in $Q_n$ has measure $p^k(1-p)^{n-k}$. If $X$ and $Y$ are upward closed, it still holds that the measure of $X\cap Y$ is at least the product of the measures of $X$ and $Y$, by the FKG inequality \cite{Gri}, or by application of a measurable map $f:Q_{bm}\to Q_m\left(\frac{u}{2^v}\right)$. Thus, we generate a counterexample in $Q_{21}$ of density $\approx 0.447$ from an approximate counterexample in $Q_7\left(\frac{3}{8}\right)$.

Consider $X$ and $Y$ to be the sets containing $1$ and $2$, respectively, in $Q_n(p)$. Let $Z$ contain the sets of size at least $l+1$, as well as the sets of size $l$ that are not in $X$ or $Y$. Then $X\cap Y^c\cap Z^c$ contains $\binom{n-2}{k-1}$ many sets of size $k$, where $1\le k\le l$. Also, $X^c\cap Y^c\cap Z$ contains $\binom{n-2}{k}$ many sets of size $k$, where $l\le k\le n-2$. The total measure of sets in exactly one of $X$, $Y$, $Z$ is thus given by
$$q=2\sum_{k=1}^{l}\binom{n-2}{k-1}p^k(1-p)^{n-k}+\sum_{k=l}^{n-2}\binom{n-2}{k}p^k(1-p)^{n-k}.$$
Compare this to the binomial expansion of $(p+(1-p))^{n-2}$:
$$\sum_{k=1}^{l}\binom{n-2}{k-1}p^{k-1}(1-p)^{n-k-1}+\sum_{k=l}^{n-2}\binom{n-2}{k}p^k(1-p)^{n-k-2}=1.$$
The first summand is $2p(1-p)$ as much, while the second summand is $(1-p)^2$ as much. Specifically, if $p=\frac13$, these coefficients both equal $\frac49$.

Since $q(p)=\frac49$ if $p=\frac13$, unless $q$ attains a local maximum at $\frac13$, a different value of $p$ could yield a larger value of $q$. For example, in the case of $n=7$ and $l=3$, $q(p)>0.447$ if $p=\frac38$. There, the measure of $Z$ is only $\approx 0.323$, but including all the supersets of $\{1,2\}$ from $X\cap Y$ boosts the measure of $Z$ up to $\approx 0.377>\frac38$.

A measure-preserving map $f:Q_{21}\to Q_7\left(\frac38\right)$ arises from the upward closed $I=\{\{1,2\},\{1,3\},\{1,2,3\}\}\subset Q_3$. Viewing $Q_{21}\cong (Q_3)^7$, we map $I$ to $1$ and $I^c$ to $0$ at every coordinate. Clearly, the preimages of $X$, $Y$, and $Z$ are upward closed.

In $Q_{21}$, the upward closed system $f^{-1}(Z)$ can be modified by including sets that map to a superset of $\{1,2\}$. Including only a maximal set in each step, the size of the set system can be incremented until it equals $\frac38$ exactly. Meanwhile, the measure of sets contained in exactly one of our three set systems is $>0.447>\frac{4}{9}$.

\section{Conclusion}

Regarding the density of $S_1$, there is a wide gap between our upper bound of $\approx 0.515$ and our counterexample of $\approx 0.447$. Nevertheless, Bhargav Narayanan, Jeff Kahn, and Sophie Spirkl \cite{Nar} also ask whether there is a counterexample where the sets $X\cap Y^c\cap Z^c$, $X^c\cap Y\cap Z^c$, and $X^c\cap Y^c\cap Z$ each have density more than $\frac{4}{27}\approx 0.148$. Moreover, can their densities each exceed $0.153$? These questions continue to be open.

\section{Acknowledgements}

The author is grateful to Imre Leader for proofreading and fruitful suggestions. The author is a recipient of the Internal Graduate Studentship in Trinity College, Cambridge, UK.

\textsc{Department of Pure Mathematics and Mathematical Statistics, University of Cambridge, Wilberforce Road, Cambridge CB3 0WB.} \\ \\
\textit{E-mail address:} \texttt{kkw25@cam.ac.uk}

\end{document}